\documentclass[reqno,11pt]{amsart}

\usepackage{color}

\newtheorem{theorem}{Theorem}[section]
\newtheorem{lemma}[theorem]{Lemma}
\newtheorem{corollary}[theorem]{Corollary}

\theoremstyle{definition}
\newtheorem{assumption}[theorem]{Assumption}

\theoremstyle{remark}
\newtheorem{remark}[theorem]{Remark}
\newtheorem{example}[theorem]{Example}

\makeatletter 
 \def\dashint{%
 \operatorname%
 {\,\,\text{\bf--}\kern-.98em\DOTSI\intop\ilimits@\!\!}}
\def\dashnorm{\,\,\text{\bf--}\kern-.5em\|}
\def\ninf{\qopname\relax\@empty{inf\phantom{p}\!\!\!}}

 \makeatother

\newcommand\sfb{{\sf b}}

\def\sft{{\sf t}}

\newcommand\bC{\mathbb{C}}

\newcommand\bR{\mathbb{R}}
\newcommand\bS{\mathbb{S}}

\newcommand\cF{\mathcal{F}}

\newcommand \cL{\mathcal{L}}

\newcommand\frN{\mathfrak{N}}

\newcommand{\loc}{{\rm loc}\,}

 \newcommand{\mysection}[1]{\section{#1}
 \setcounter{equation}{0}}

\newcommand{\nliminf}{\operatornamewithlimits{\underline{lim}}}

\begin{document}

\title[On weak solutions]
{Once again on weak solutions of time inhomogeneous
It\^o's equations with VMO diffusion and Morrey drift}
\author{N.V. Krylov}

\email{nkrylov@umn.edu}
\address{School of Mathematics, University of Minnesota, Minneapolis, MN, 55455}
 
 

\begin{abstract} 
We prove the existence and weak uniqueness of weak solutions
of It\^o's stochastic time dependent equations with
irregular diffusion and drift terms
of Morrey class with mixed norms.  
\end{abstract}

\maketitle

\mysection{Introduction}
                                                  \label{section 3.11.1}

This paper is a natural complement of \cite{Kr_pt1}
where the drift term was assumed to be the sum
of two terms, one of which was bounded in $x$
and $L_{2}$ in $t$, and another one was in a Morrey
class for each $t$ with small norm. In this paper
we concentrate on the case when the drift is
in a Morrey class with respect to $(t,x)$
with mixed norms. Our results substantially extend
those in \cite{Kr_22-2}.

Let $\bR^{d}$ be a $d-$dimensional Euclidean space of points
$x=(x^{1},...,x^{d})$ with $d\geq2$. Let $(\Omega,\cF,P)$ be a 
complete probability space,
carrying a $d $-dimensional Wiener process
  $w_{t}$. Fix $\delta\in(0,1]$ and denote by $\bS_{\delta}$ the set of $d\times d$ symmetric matrices
whose eigenvalues lie in $[\delta,\delta^{-1}]$.

Assume that on $\bR^{d+1}=\{(t,x):t\in\bR,x\in\bR^{d}\}$ we are given Borel $\bR^{d}$-valued 
function  $b=(b^{i})$ and $\bS_{\delta}$-valued $\sigma 
=(\sigma^{ij})$.
We are going
to   investigate
the equation
\begin{equation}
                                         \label{6.15.2}
 x_{s}=x +\int_{0}^{s}\sigma (t+u,x_{u})\,dw _{u}+
\int_{0}^{s}b(t+u,x_{u})\,du.
\end{equation}

We are interested in the so-called weak solutions, that is
solutions that are not necessarily $\cF^{w}_{s}$-measurable,
where $\cF^{w}_{s}$ is the completion of $\sigma(w_{u}:u\leq s)$. We present
sufficient conditions for the equation to have such solutions on appropriate probability spaces and investigate uniqueness of their distributions.

We just reproduced part of the introduction 
in \cite{Kr_pt1}. The reader interested in learning
more about the history of the problem, motivation, and the literature
is sent to 
   \cite{BFGM_19}, \cite{Ki_23}, \cite{RZ_20}     and to the
the introduction in \cite{Kr_pt1}.

Our Morrey type condition on $b$ is stated
in terms of mixed norms (different powers of summability with respect $t$ and $x$).
In \cite{Ki_23} the Morrey type condition 
on $b$ is stated in terms of  $L_{p}$-norms
in $(t,x)$ and the weak solvability of \eqref{6.15.2}
is proved, if $\sigma$ is the unit matrix, along with weak uniqueness provided
the solutions possess some additional properties.
Our result in this respect contains \cite{Ki_23}
and also allow us to give conditions
for unconditional weak uniqueness.

  Still our general uniqueness theorem and uniqueness
theorems in \cite{RZ_20} are   conditional.
We prove uniqueness only in the class of solutions
(which is proved to be nonempty) admitting certain
estimates, however, as we said before, there are cases in which we prove unconditional weak uniqueness.

Our $\sigma$ is not constant or continuous
and 
it is worth saying
that
restricting
the situation to the one when $\sigma$ and $b$
are independent of time allows one to
relax the  conditions on $b$ significantly
further, see, for instance, \cite{KS_19}
and the 
references therein.

In Remark 
\ref{remark 10.31.1}  we compare our results with some
of those in an excellent papers by
R\"ockner and Zhao \cite{RZ_20} and \cite{Ki_23}. 
  By the way,
G. Zhao (\cite{Zh_20_1}) gave an example showing
that, if in the definition of $\mathfrak{b}_{\rho}$   we replace $r$ with $r^{ \alpha}$, $\alpha>1$, 
and keep \eqref{9.27.4} (with $k=0$),  weak uniqueness
may fail even in the time homogeneous case
and unit diffusion. In Example \ref{example 12.21.4}
we show that in the time inhomogeneous case even
the existence may fail.

Here is an example 
 in which we prove existence and (unconditional)
weak
 uniqueness of weak solutions:
$|b|=cf$, where 
the constant $c>0$ is small enough and
$$
f(t,x)=I_{1>t>0,|x|<1}|x|^{-1}\Big(\frac{|x|}{\sqrt t}\Big)^{1/(d+1)},\,\,
\sigma=2(\delta^{ij})+I_{x \ne0}\zeta(x )\sin(\ln|\ln |x |),
$$
where $\zeta$ is any smooth symmetric $d\times d$-matrix valued 
function vanishing for $|x|>1/2$ and
satisfying $|\zeta|\leq 1$.
This example is inadmissible in  
  \cite{RZ_20} because $b$ is too singular and
  $\sigma$ is not constant but admissible in
\cite{Ki_23} if $\sigma$ is constant
and yields weak solutions conditionally weakly unique. This example
does not fit into the scheme in \cite{Kr_pt1}
(barely misses)
  because of special $b$.

 The paper is organized as follows. 
In Section \ref{section 12.24.1}
 we prove the solvability
of \eqref{6.15.2} when the drift is the sum
of terms with different summability properties.
In Section \ref{section 3.7.1} we deal with  
weak uniqueness   and construct the corresponding Markov
processes. This time the drift is not split.  
Section \ref{section 3.10.1} contains
a result from \cite{Kr_23_1} used in Section \ref{section 3.7.1}.

We conclude the introduction by some notation.
We set  
$$
 D_{i}  =\frac{\partial}{\partial x^{i}} ,\quad Du= (D_{i}u),\quad
  D_{ij} = D_{i}D_{j}  ,\quad D^{2}u=(D_{ij}u),
\quad
 \partial_{t} =\frac{\partial}{\partial t} .
$$
 If $\sigma =(\sigma^{i...})$ by $|\sigma|^{2}$
we mean the sum of squares of all entries.

 Introduce
$$
B_{R}(x)=\{y\in\bR^{d}:|x-y|<R\},\quad B_{R} =B_{R}(0),
$$
$$
C_{\tau,\rho}(t,x)=[t,t+\tau)\times B_{\rho}(x),\quad C_{\rho}...=C_{\rho^{2},\rho}...,\quad C_{\rho}=C_{\rho}(0,0),
$$
and let $\bC_{\rho}$ be the collection of
$C_{\rho}(t,x)$.

In the proofs of our results
we use various (finite) constants called $N$ which
may change from one occurrence to another
and depend on the data only in the same way
as it is
  indicated in the statements
of the results.

\mysection{Solvability of It\^o's equations}
                          \label{section 12.24.1}

Let $  d\geq 2$ and let $(\Omega,\cF,P)$ be a complete probability
space. Let $\cF_{t}, t\geq0$, be an increasing family of
complete $\sigma$-fields $\cF_{t}\subset\cF$,  
and let $w_{t}$ be an $\bR^{d}$-valued Wiener process
relative to $\cF_{t}$.

It is well known that,
if  $\sigma$ and $b$ are smooth
and $b$ is bounded, the solutions of the system
\begin{equation}
                                        \label{11.29.20}
x _{s}=x  +\int_{0}^{s}\sigma (\sft_{r},x_{r})\,dw_{r}
+\int_{0}^{s}b (\sft_{r},x_{r}) \,dr,\quad \sft_{s}=t+s
\end{equation}
form a strong Markov process $X$ with trajectories $(\sft_{s},x_{s})$. 

Define
$$
  \sfb_{\rho}=\sup_{r\leq\rho}r^{-1}
\sup_{(t,x)\in\bR^{d+1}}\sup_{C\in\bC_{r}}
E_{t,x}\int_{0}^{\tau_{C}}|b(\sft_{s},x_{s})|\,ds,
$$
where $\tau_{C}$ is the first exit time of 
$(\sft_{s},x_{s})$ from $C$.

To continue we need some notation which 
are somewhat different from what we use
in Sections \ref{section 3.7.1} and \ref{section 3.10.1}.
For $p,q\in[1,\infty)$ and domain $Q\subset\bR^{d+1}$ by $L_{p,q}(Q)$
we mean the space of Borel (real-, vector- or matrix-valued)
 functions on $Q$   with finite norm given by
$$
\|f\|_{L_{p,q}(Q)}^{q}=\|fI_{Q}\|_{L_{p,q}}^{q}
=\int_{\bR}\Big(\int_{\bR^{d}}|fI_{Q}(t,x)|^{p}\,
dx\Big)^{q/p}\,dt
$$
if $p\geq q$ and by
$$
\|f\|_{L_{p,q}(Q)}^{p}=\|fI_{Q}\|_{L_{p,q}}^{p}
=\int_{\bR^{d}}\Big(\int_{\bR}|fI_{Q}(t,x)|^{q}\,dt
 \Big)^{p/q}\,dx
$$
if $p\leq q$. 
Set $L_{p,q}=L_{p,q}(\bR^{d+1})$.
These definitions extend naturally when one or both $p,q$ are infinite.
  As usual, we write something like $f\in L_{p,q,\loc}$
if $f\zeta\in L_{p,q}$ for any 
infinitely differentiable $\zeta$ with compact support. We write $\|u,v,...\|_{L_{p,q}}$ to mean
the sum of the $L_{p,q}$-norms of what is inside.

By $W^{1,2}_{p,q}(Q)$ we mean the collection
of $u$ such that $\partial_{t}u$, $D^{2}u $, $Du$, $u
\in L_{p,q}(Q)$. The norm in $W^{1,2}_{p}(Q)$
is introduced in an obvious way.
We abbreviate $W^{1,2}_{p,q } =W^{1,2}_{p,q}(\bR^{d+1})$.

If a Borel $\Gamma\subset \bR^{d+1}$, by $|\Gamma|$ we mean its Lebesgue
measure and
$$
\dashint_{\Gamma}f(x)\,dxdt=\frac{1}{|\Gamma|}
\int_{\Gamma}f(x)\,dxdt.
$$
 If $C\in\bC_{\rho}$ we set
$$
\dashnorm f\|_{L_{p,q}(C)}=\|1\|_{L_{p,q}(C)}^{-1}
\|f\|_{L_{p,q}(C)}=N(d)\rho^{-d/p-2/q}\|f\|_{L_{p,q}(C)}.
$$

Take the Fabes-Stroock constant $d_{0}=d_{0}(d,\delta)\in(d/2,d)$
 introduced in \cite{Kr_pta} and let us
say that $(p,q)$ are {\em admissible\/} if 
$$
p ,q  \in[1,\infty],\quad
 \frac{d_{0}}{p }+\frac{1}{q }\leq1.
$$
Also take $m_{b}=m_{b}(d,\delta)>0$
introduced in \cite{Kr_pta}.

Here is a generalization of Corollary 2.14 
of \cite{Kr_21_1}.

\begin{theorem}
                     \label{theorem 9.27.10} 
Assume that  $\sigma$ and $b$ are smooth
and $b$ is bounded and there is a nonnegative integer $k$ and
there are Borel functions $b_{i}(t,x)$,
$0\leq i\leq k$, 
such that $b =\sum_{i=0}^{k}b_{i} $,  and 
we are given admissible $(p_{i},q_{i})$,
$i\leq k$. Define  
$$
\mathfrak{b}_{\rho}=\sup_{r\leq\rho}r
\sup_{C\in \bC_{r}}\sum_{i=0}^{k}
\dashnorm b_{i}\|_{L_{p_{i},q_{i}}(C)},
$$
introduce  
$\hat b=
\hat b(d,\delta)$ so that $ N\hat b= m_{b}/4$, where $  N=N(d,\delta)$ is taken from   Theorem 1.1 
of \cite{Kr_21_1} and suppose that  
\begin{equation}
                 \label{9.27.4}   
\mathfrak{b}_{\rho_{b}}\leq \hat b   
\end{equation}
holds
for some $\rho_{b}\in(0,\infty)$. Then
\begin{equation} 
                          \label{9.27.5}
 \sfb_{ \rho_{b}/2}\leq m_{b}.
\end{equation}
\end{theorem}

One proves this theorem by   repeating
the proof of Theorem 1.1 of \cite{Kr_21_1},
where the fact that there $k=0$ was not used at all,
and also using the argument in Step 1 of the proof
of Theorem 1.2 of \cite{Kr_21_1}.
After that the same argument as in
Corollary 2.14 
of \cite{Kr_21_1} yields the result.

\begin{example}
                \label{example 1.27.1}
One of situations when $\mathfrak{b}_{\rho}$ is finite presents when
$k=1$,
$|b_{0}(t,x)|\leq h_{0}(x)$, $|b_{1}(t,x)|
\leq h_{1}(t)$ and, say $h_{0}(x)\leq c |x|^{-1}$,
where $c$ is sufficiently small, and $h_{1}\in L_{2}(\bR)$.
In that case one can take $p_{1}=d_{0},
q_{1}=\infty$, $p_{2}=\infty,q_{2}=1$.

Indeed, if $|x_{0}|\leq 2r$, then
$$
\dashint_{B_{r}(x_{0})}|x|^{-d_{0}}\,dx
\leq 2^{d}\dashint_{B_{2r} }|x|^{-d_{0}}\,dx=N(d)r^{-d_{0}}.
$$
and if $|x_{0}|\geq 2r$, then 
$|x|^{-1}\leq r^{-1}$ on $B_{r}(x_{0})$
and $\dashnorm |\cdot|^{-1}\|_{L_{d_{0}}
(B_{r}(x_{0}))}\leq r^{-1}$.

Also    
$$
\dashint_{s}^{s+r^{2}}h_{1}(t)\,dt
\leq r^{-1}\Big(\int_{s}^{s+r^{2}}
h_{1}^{2}(t)\,dt\Big)^{1/2}
$$
and the integral here tends to zero
as $r\downarrow 0$ uniformly with respect to $s$.
Therefore, by taking $c$ small enough
 and taking appropriately small $\rho_{b}$ we can satisfy
\eqref{9.27.4} with any given $\hat b>0$.
\end{example}
 
In the following theorem we prove
the existence of weak solutions of
equation \eqref{11.29.20}. Somewhat
unusual split in its assumption about $b$
is caused by the necessity to use
smooth approximations of the $b_{i}$'s
converging to the $b_{i}$'s in the corresponding norms.
\begin{theorem}
               \label{theorem 10.21.1} 
Suppose that  
\eqref{9.27.4} holds for some $\rho_{b}\in(0,\infty)$   with admissible
$p_{i},q_{i} $  such that,
for each $i$, either (a) 
$p_{i}+q_{i}<\infty$, or (b) $p_{i}<\infty$,
$q_{i}=\infty$ and $b_{i}$ is independent of $t$.
Then 

(i) there is a probability space 
$(\Omega ,\cF ,P )$,
a filtration of $\sigma$-fields $\cF _{s}\subset \cF $, $s\geq0$,
a process $w _{s}$, $s\geq0$, which is a $d$-dimensional Wiener process
relative to $\{\cF _{s}\}$, and an $\cF _{s}$-adapted
process $x_{s}$ such that 
 (a.s.) for all   $s\geq0$ equation \eqref{11.29.20} holds with $(t,x)=(0,0)$.

(ii) Furthermore, for any 
  nonnegative Borel $g$ on $\bR^{d}$ and $f$ on $\bR^{d+1} $ and
 $ T\in(0,\infty)$   we have
\begin{equation}
                                    \label{3.8.1}
E   \int_{0}^{T}  
f(s,x _{s} )\,ds \leq N(d,\delta, T,\rho_{b})  
 \|  f\| _{L_{p_{i},q_{i}}} ,
\end{equation}
\begin{equation}   
                    \label{3.8.2} 
E   \int_{0}^{T}  
g(x _{s} )\,ds \leq N(d,\delta, T,\rho_{b})  
 \|  g\| _{L_{d_{0}}(\bR^{d})} .
\end{equation}

\end{theorem} 

Proof. Approximate $\sigma,b$ by smooth
$\sigma^{(\varepsilon)},b^{(\varepsilon)}$,
by using mollifying kernel $\varepsilon^{-d-1}\zeta
(t/\varepsilon,x/\varepsilon)$, where nonnegative
$\zeta\in C^{\infty}_{0}(\bR^{d+1})$ has unit integral and $\zeta(0)=1$. Then set
$b_{i}^{\varepsilon}(t,x)=b_{i}^{(\varepsilon)}(t,x)\zeta(
\varepsilon t,\varepsilon x)$ to make the new $b_{i}$
have compact support. Observe that $b_{i}^{\varepsilon}$
satisfy   \eqref{9.27.4}
with the same $\hat b,\rho_{b} $.
Therefore, the corresponding Markov process
$(\sft_{t},x^{\varepsilon}_{t})$
satisfies $ \sfb^{\varepsilon}_{\rho_{b}/2}\leq m_{b}$
which makes available 
all results of \cite{Kr_pta}.
In particular, by Corollary 3.10 of \cite{Kr_pta}
for any $\varepsilon,n>0$ and
  $r>s\geq 0$  
\begin{equation}
                                    \label{11.24.4}
E_{0,0}\sup_{u\in[s,r]}|x^{\varepsilon}_{u}-x^{\varepsilon}_{s}|^{ n}
\leq N(  |r-s| ^{ n/2}+|r-s| ^{ n}),
\end{equation}
where $N=N(n, \rho_{b},d,\delta)$. This implies that the  $P_{0,0}$-distributions of $x^{\varepsilon}_{\cdot}$
are precompact on $C([0,\infty),\bR^{d})$
and a subsequence as $\varepsilon=\varepsilon_{n}
\downarrow 0$ of them converges
to the distribution of a process $x^{0}_{\cdot}$ defined
on a probability space (the coordinate process on $\Omega=C([0,\infty),\bR^{d})$ 
with cylindrical $\sigma$-field $\cF$ 
completed with respect to
 $P$, that is
 the limiting distribution of $x^{\varepsilon}_{\cdot}$). Furthermore, by Theorem 5.9 (ii) of \cite{Kr_pta} for any 
  nonnegative Borel $g$ on $\bR^{d}$ and $f$ on $\bR^{d+1} $ and
 $\varepsilon,T\in(0,\infty)$   we have
\begin{equation}
                                    \label{11.24.3} 
E _{0,0} \int_{0}^{T}  
f(s,x^{\varepsilon}_{s} )\,ds \leq N(d,\delta, T,\rho_{b})  
 \|  f\| _{L_{p_{i},q_{i}}} ,
\end{equation}
\begin{equation}   
                    \label{1.27.4}
E _{0,0} \int_{0}^{T}  
g(x^{\varepsilon}_{s} )\,ds \leq N(d,\delta, T,\rho_{b})  
 \|  g\| _{L_{d_{0}}(\bR^{d})} ,
\end{equation}
which by continuity is extended to $\varepsilon=0$
for bounded continuos $f$ and then by the usual measure-theoretic argument for all Borel $f\geq0$.   
This proves (ii).

 After that arguing as in the proof of Theorem 3.9
of \cite{Kr_pt1}  proves assertion (i). Here passing to the limit
in the drift term the case (a) we use \eqref{11.24.3} and in the case (b)
we use \eqref{1.27.4}.
The theorem is proved. \qed  

\begin{remark}
                      \label{remark 3.11.1}
Actually as is easy to see, in case (b) the condition that $b_{i}$
is independent of $t$ can be replaced with
the following which is somewhat cumbersome: 
for any $R\in(0,\infty)$
$$
\lim_{\varepsilon\downarrow 0}\int_{B_{R}}\sup_{t}
|b^{(\varepsilon)}-b|^{d_{0}}(t,x)\,dx=0.
$$
\end{remark}

\begin{remark}
                 \label{remark 1.26.1}
It may look like assertion (i)
of Theorem \ref{theorem 10.21.1}
is a generalization of Theorem 3.1 (i)
of \cite{Kr_20} about the solvability of
\eqref{11.29.20} with $b\in L_{p,q}$
and $d/p+1/q\leq 1$. However,
in the typical case of $k=0$,
along with $b \in L_{p_{0} ,q_{0} ,\loc}$,
$d_{0}/p_{0} +1/q _{0}\leq1$,
we require \eqref{9.27.4} to hold
and, if we ask ourselves what
$p,q$ should be in order the
inclusion $b\in L_{p,q}$ to imply 
\eqref{9.27.4}, the answer is 
$d/p+2/q\leq 1$, somewhat disappointing.
At the same time in the next example 
we show that Theorem 3.1 (i)
of \cite{Kr_20} does not cover
all applications of Theorem \ref{theorem 10.21.1}.
\end{remark}

In assumption  
\eqref{9.27.4}   the size
of $\hat b$ could not be too large.

\begin{example}
                           \label{example 12.21.01}
Let
$$
b(t,x)=b(x)=-\frac{d}{|x|}\frac{x}{|x|}I_{x\ne0}, \quad \sigma=\sqrt2 (\delta^{ij}).
$$
Then as is easy to see, for any $p\in(d_{0},d)$
and any $q$ the quantity $\rho\dashnorm b\|_{L_{p,q}(C)}$, $\rho>0,C\in\bC_{\rho}$, is bounded.
However, the equation $dx_{t}=\sigma\,dw_{t}
+b(x_{t})\,dt$ with initial condition $x_{0}=0$
does not have any solution.

Indeed, if it does, then by It\^o's formula
\begin{equation}
                              \label{12.21.5}
|x_{t}|^{2}=2d\int_{0}^{t}I_{x_{s}=0}\,ds
+2\sqrt2\int_{0}^{t}x_{t}\,dw_{t}.
\end{equation}
Here the first integral is the time spent at the
origin by $x_{s}$ up to time $t$. This integral is zero, because by using It\^o's formula for
$|x^{1}_{t}|$, one sees that the local time
of $x^{1}_{t}$ at zero exists and is finite,
implying that the real time spent at zero is zero.

Then \eqref{12.21.5} says that the local martingale
starting at zero which stands on the right is
nonnegative. But then it is identically zero,
implying the same for $x_{t}$. However,
$x_{t}\equiv0$, obviously, does not satisfy our equation.

At the same time according to Theorem
\ref{theorem 10.21.1},
the equation $dx_{t}=\sigma\,dw_{t}
+\varepsilon b(x_{t})\,dt$ with initial condition $x_{0}=0$ does have solutions if $\varepsilon$ is sufficiently small.
Observe that $b\not \in L_{p,q,\loc}$
for any $p,q\in(0,\infty)$ satisfying
$d/p+1/q\leq1$, so this example is not covered by Theorem 3.1 (i)
of \cite{Kr_20}.
\end{example}

It turns out that in the definition of $\mathfrak{b}$ one
cannot replace $r$ with $r^{1+\alpha}$,
no matter how small $\alpha>0$ is.

\begin{example}
                         \label{example 12.21.4}

Take   numbers $\alpha$ and $\beta$ satisfying
$$
0<\alpha\leq \beta <1,\quad \alpha+\beta=1
$$
 and set
$$
b(t,x)=-\frac{1}{t^{\alpha}|x |^{\beta}}\frac{x }{|x |}
I_{0<|x|\leq 1,t\leq 1}.
$$
Using that $d_{0}<d$, it is not hard to find $p,q$
such that  $d_{0}/p+1/q<1$ and the quantity
$\rho^{1+\alpha}\dashnorm b\|_{L_{p,q}(C)}$,
$\rho>0, C\in\bC_{\rho}$, is bounded.
However, as we know from \cite{Kr_20}, the equation $dx_{t}=dw_{t}+\varepsilon
b(t,x_{t})\,dt$ with zero initial condition
does not have solutions no matter how small $\varepsilon>0$ is (actually $\varepsilon=1$ in \cite{Kr_20} but self-similar transformations take care of any $\varepsilon>0$).
\end{example}

\begin{remark}
                                \label{remark 3.9.1}
If $b\equiv0$, it turns out
that for any admissible $(p,q)$, $R\in(0,\infty)$,
$x\in \bR^{d}$ and Borel $f(t,x)\geq 0$
\begin{equation}
                             \label{3.9.2}
E\int_{0}^{\tau}f(s,x_{s})\,ds
\leq N(d,\delta)R^{2}\dashnorm f\|_{L_{p,q}(C_{R}(0,x))},
\end{equation}
where $\tau$ is the first exit time of $(s,x_{s})$
from $C_{R}(0,x)$.

Indeed, if $R=1$, this follows from \eqref{3.8.1}, where we take $T=1$,
any appropriate $\rho_{b}$ and observe that
$\tau\leq 1$ and we may assume that $f=0$
outside $C_{1}(0,x)$. The case of general $R$
is treated by parabolic scaling of $\bR^{d+1}$.
\end{remark}

This simple observation has the following implication
in which
$$
\cL_{0} u(t,x)=\partial_{t}u+a^{ij}(t,x)D_{ij}u(t,x),\quad a=(1/2)\sigma^{2}.
$$
\begin{lemma}
                          \label{lemma 3.9.2}
Let $(p,q)$ be admissible and finite, $x\in\bR^{d}$,
$R\in(0,\infty)$, $u\in W^{1,2}_{p,q}(C_{R}(0,x))$
and $u=0 $ on $\partial'C_{R}(0,x)$ (that is
$\big(\partial B_{R}(x)\times[0,R^{2}]\big)
\cup\big(\bar B_{R}(x)\times \{R^{2}\}\big)$).
Then
\begin{equation}
                             \label{3.9.3}
|u(0,0)| 
\leq N(d,\delta)R^{2}\dashnorm \cL_{0}u\|_{L_{p,q}(C_{R}(0,x))}.
\end{equation}
\end{lemma}

Proof. First note that, since $d_{0}>d/2$ we have
$d/p+2/q<2$ and $u$ is continuous in 
$\bar C_{R}(0,x)$ by embedding theorems. Then
approximate $u$ in $W^{1,2}_{p,q}$-norm by smooth
functions $u^{n}$ vanishing on $\partial'C_{R}(0,x)$.
By It\^o's formula
$$
u^{n}(0,0)=-E\int_{0}^{\tau}\cL_{0}u^{n}
(s,x_{s})\,ds.
$$
In light of \eqref{3.9.2} estimate \eqref{3.9.3}
holds with $u^{n}$ in place of $u$. Sending
$n\to\infty$ yields \eqref{3.9.3} as is and
proves the lemma. \qed

Here is It\^o's formula we have on the basis
of Theorem \ref{theorem 10.21.1}. 

\begin{theorem}
                      \label{theorem 12.11.3}
(i) Suppose that $k=0$ and
\eqref{9.27.4} holds for some $\rho_{b}\in(0,\infty)$ 
and  
\begin{equation}
                               \label{12.13.4}
p_{0},q_{0}\in (1,\infty),\quad
\frac{1}{2}\leq\frac{1}{\beta_{0}}:=\frac{d_{0}}{p_{0}}+\frac{1}{q_{0}}<1. 
\end{equation}

(ii) Let $x_{t}$ be a solution of equation \eqref{11.29.20}   with $(t,x)=(0,0)$ on a probability space such that

(a)
for $p=  p_{0}/\beta_{0}$,
$q=  q_{0}/\beta_{0}$, any $R\in(0,\infty)$
and Borel nonnegative $f$ on $\bR^{d+1}$
\begin{equation}
                             \label{12.11.6}
E\int_{0}^{\tau_{R}}f(s,x_{s})\,ds\leq N\|f\|_{L_{p,q}},
\end{equation}
where $N$ is independent of $f$ and $\tau_{R}$ is the first exit time of $(s,x_{s})$ from $C_{R}$.  

(iii) Let
  $u\in W^{1,2}_{p,q}(C_{R})$ be such that 
  $Du\in L_{r,k}(C_{R})$, where $(r,k)=
(\beta_{0}-1)^{-1}(p_{0},q_{0})$.

Then, with probability one for all $t\geq0$,
$$
u(t\wedge\tau_{R},x_{t\wedge\tau_{R}})
=u( 0)+\int_{0}^{t\wedge\tau_{R}}D_{i}u 
\sigma^{ik} (s,x_{s})\,dw^{k}_{s}
$$
$$
+\int_{0}^{t\wedge\tau_{R}}[
\partial_{t}u(s,x_{s})+ a^{ij} D_{ij}u(s,x_{s})
+b^{i} D_{i}u(s,x_{s})]\,ds
$$
and the stochastic integral above is a square-integrable
martingale, where $\tau_{R}$ is the first
exit time of $x_{t}$ from $B_{R}$.
\end{theorem}

Proof. The last statement, of course,
follows from \eqref{12.11.6} 
and the fact that $2p\leq r,2q\leq k$.
To prove the rest we   approximate  $u$
by smooth functions $u^{(\varepsilon)}=\zeta_{\varepsilon}*u$, where $\zeta_{\varepsilon}
(t,x)=\varepsilon^{-d-2}\zeta(t/\varepsilon^{2},
x/\varepsilon)$, and $\zeta\in C^{\infty}_{0}
(\bR^{d+1})$ has support in $(-1,0)\times B_{1}$
and unit integral. Since $d/p+2/q<2$ ($d<2d_{0}$),
by embedding theorems $u\in C(\bar C_{R})$ and, therefore, $u^{(\varepsilon)}\to u$ as $\varepsilon\downarrow 0$ uniformly
in any $C_{R'}$ with $R'<R$.

Fix $R'<R$. Then for all sufficiently small
$\varepsilon>0$  by It\^o's formula
$$
u^{(\varepsilon)}(t\wedge\tau_{R'},x_{t\wedge\tau_{R'}})
=u^{(\varepsilon)}( 0)+\int_{0}^{t\wedge\tau_{R'}}D_{i}u^{(\varepsilon)} 
\sigma^{ik} (s,x_{s})\,dw^{k}_{s}
$$
\begin{equation}
                                      \label{10.31.1}
+\int_{0}^{t\wedge\tau_{R'}}[
\partial_{t}u^{(\varepsilon)}+ a^{ij} D_{ij}u^{(\varepsilon)}
+b^{i} D_{i}u^{(\varepsilon)} ](s,x_{s})\,ds.
\end{equation}
We send $\varepsilon\downarrow 0$ and observe that
$u^{(\varepsilon)}\to u$ in $W^{1,2}_{p,q}(C_{R'})$
and $Du^{(\varepsilon)}\to Du$ in $L_{r,k}(C_{R'})$. Hence,
\eqref{12.11.6} allows us easily to pass
to the limit in \eqref{10.31.1}, for instance,
by using H\"older's inequality   we obtain 
\begin{equation}
                               \label{12.24.3}
\|gh\|_{L_{p,q}(C_{\rho})}
\leq \|g\|_{L_{p_{0},q_{0}}(C_{\rho})}
\|h\|_{L_{r,s}(C_{\rho})},
\end{equation}
implying that
$$
E\int_{0}^{\tau_{R'}}|b|\,|Du-Du^{(\varepsilon)}|
(t,x_{t})\,dt\leq N\||b|\,|Du-Du^{(\varepsilon)}|\|_{L_{p,q}(C_{R'})}
$$
$$
\leq N\|b\|_{L_{p_{0},q_{0}}(C_{R})}
\|Du-Du^{(\varepsilon)}\|_{L_{r,s}(C_{R'})}\to0.
$$
It follows that \eqref{10.31.1} holds with
$u$ in place of $u^{(\varepsilon)}$. After that 
it only remains to send $R'\uparrow R$ and again use 
\eqref{12.11.6}. The theorem is proved. \qed

\begin{remark}
                        \label{remark 11.28.1}
The assumption that $Du\in L_{r,s}(C_{R})$
looks unrealistic because Sobolev embedding
theorems  do not provide such high integrability 
of $Du$ for functions $u\in W^{1,2}_{p,q}(C_{R})$. However, if $u$ is in the Morrey
class $E^{1,2}_{p,q,\beta}(C_{R})$, then
$Du\in L_{r,k}(C_{R})$  indeed (cf. Remark \ref{remark 10.27.2}).
\end{remark}

\mysection{Weak uniqueness  
and a Markov process}
                             \label{section 3.7.1}

Here we prove a generalization of the Stroock-Varadhan theorem in \cite{SV_79} obtained
for $\sigma$ which is uniformly continuous in $x$
uniformly in $t$
and bounded $b$.
  We need an additional assumption on $a$
and can relax conditions imposed on $b$ in Section
\ref{section 12.24.1}.
Since $a$ will have some regularity the range of $p_{0},q_{0}$ can be substantially extended.
Indeed, observe that if $ d_{0}/p_{0}
+1/q_{0}= 1$, then $ 1<d/p_{0}+2/q_{0}< 2$
since $d>d_{0}>d/2 $ (cf.~\eqref{3.14.1}).

An important distinction of the rest of the article
from Section \ref{section 12.24.1} is that
here (and in Section \ref{section 3.10.1})
for $p,q\in[1,\infty)$ and domain $Q\subset\bR^{d+1}$ by $L_{p,q}(Q)$
we mean the space of Borel (real-, vector- or matrix-valued)
 functions on $Q$   with finite norm given in
{\em one of two ways which is fixed throughout the 
rest of the paper\/}:
\begin{equation}
                                 \label{3.27.3}
\|f\|_{L_{p,q}(Q)}^{q}=\|fI_{Q}\|_{L_{p,q}}^{q}
=\int_{\bR}\Big(\int_{\bR^{d}}|fI_{Q}(t,x)|^{p}\,
dx\Big)^{q/p}\,dt
\end{equation}
or
\begin{equation}
                                 \label{4.3.2}
\|f\|_{L_{p,q}(Q)}^{p}=\|fI_{Q}\|_{L_{p,q}}^{p}
=\int_{\bR^{d}}\Big(\int_{\bR}|fI_{Q}(t,x)|^{q}\,dt
 \Big)^{p/q}\,dx.
\end{equation}

Naturally,  $\dashnorm\cdot\|_{L_{p,q}(Q)}$
and the spaces $W^{1,2}_{p,q}(Q)$ are now introduced
in the same way as in Section \ref{section 12.24.1} but with the new
meaning of $L_{p,q}(Q)$.

Fix $p_{0},q_{0},\beta_{0},\beta_{0}'$ such that
\begin{equation}
                           \label{3.14.1}
\beta_{0}\in(1,2),\quad
\beta_{0}'\in(1,\beta_{0}),\quad p_{0},q_{0}\in(\beta_{0}',\infty),\quad 
\frac{d}{p_{0}}
+\frac{2}{q_{0}}\geq 1, \quad.
\end{equation}

Take an $\alpha\in(0,1)$ and $\theta(d,\delta,p,q,\alpha)$ introduced in Assumption
\ref{assumption 3.25.1}  and set
$$
\tilde \theta=\theta(d,\delta,p_{0}/\beta_{0},q_{0}/\beta_{0},\alpha)\wedge
\theta(d,\delta,p_{0}/\beta'_{0},q_{0}/\beta'_{0},\alpha)
$$

Also take some $\rho_{a},\rho_{b}\in(0,1] $, take
$\check b (d,\delta, p, 
q, \rho_{a}, \beta_{0},\alpha )$  from Theorem \ref{theorem 10.27.1} and set  
$$
\tilde b=
\check b(d,\delta, p_{0}/\beta_{0}, 
q_{0}/\beta_{0}, \rho_{a}, \beta_{0},\alpha )\wedge
\check b(d,\delta, p_{0}/\beta'_{0}, 
q_{0}/\beta'_{0}, \rho_{a} ,\beta_{0}',\alpha ).
$$
Finally,  one more restriction 
on the drift term is related to the following
condition:
\begin{equation}
                         \label{4.15.3}
 \frac{d}{p_{0}}+\frac{1}{q_{0}}\leq 1\,\,\text{and}\,\,
\begin{cases}\text{either}&p_{0}\geq q_{0}\,\,
\text{and}\,\, L_{p,q}\,\,
\text{is defined as in} \,\eqref{3.27.3},\\
\text{or}&p_{0}\leq q_{0}\,\,\text{and}\,\, L_{p,q}\,\,
\text{is defined as in} \,\eqref{4.3.2}.
\end{cases}
\end{equation}
 Generally regardless of condition 
\eqref{4.15.3}, define $\hat b\leq 1$ from $N_{1} \hat b\leq m_{b}$ 
 or, if condition 
\eqref{4.15.3} is satisfied, from $
  N  \hat b\leq m_{b}$,
where $N=N (d,\delta,p_{0},q_{0})$ is taken
from \eqref{3.8.3} and $N_{1} $
depending only on
$d,\delta,p_{0},q_{0},\beta_{0},\rho_{a},\alpha$
is taken from \eqref{3.14.8}. 

Throughout this  section we suppose that the following assumptions are satisfied.

\begin{assumption}
                        \label{assumption 12.12.3}
We have
$$
 a^{\sharp}_{x,\rho_{a}}:=\sup_{\substack{\rho\leq\rho_{a}\\C\in\bC_{\rho}}}\dashint_{C}|a (t,x)-a_{C}(t)|
\,dx dt \leq \tilde \theta,
$$
where
$$
a_{C}(t)=\dashint_{C}a (t,x)\,dxds\quad (\text{note}\quad t\quad \text{and}\quad ds).
$$
\end{assumption}

\begin{assumption}
                      \label{assumption 3.14.1} 
We have
\begin{equation}
                           \label{3.14.2}
 \mathfrak{b}_{\rho_{b}}\leq \hat b\wedge \tilde b ,
\end{equation}
where
$$
\mathfrak{b}_{\rho}=\sup_{r\leq\rho}r
\sup_{C\in \bC_{r}} 
\dashnorm b \|_{L_{p_{0},q_{0}}(C)}.
$$

\end{assumption}

The following is very important.

\begin{remark}
                            \label{remark 4.2.1}
Consider equation \eqref{11.29.20} with zero
initial data and make the change of variables
$ x_{t} = \rho _{b}y_{\rho_{b}^{-2}t}  $, 
$B_{t}=\rho_{b} w_{\rho_{b}^{-2}t}$. Then
\begin{equation}
                                    \label{4.4.1}
dy_{t} = \tilde  b(t,y_{t} )\,dt+\tilde 
\sigma( t,y_{t} )\,dB_{t},
\end{equation}
where $\tilde b(t,x)=\rho _{b}b(\rho_{b}^{ 2}t,\rho_{b}x)$, $\tilde\sigma(t,x)
=\sigma(\rho_{b}^{ 2}t,\rho_{b}x)$, and $B_{t}$
is a Wiener process.

Taking into account that $\rho_{b}\leq1$, it is easy to check that $\tilde \sigma$
and $\tilde b$ satisfy Assumptions \ref{assumption 12.12.3}
and \ref{assumption 3.14.1} with the {\em same\/}
$\rho_{a},\tilde \theta,\hat b,\tilde b$ and 1
in place of $\rho_{b}$ in \eqref{3.14.2}.
At the same time the issues of existence and uniqueness of solutions of \eqref{4.4.1} and
\eqref{11.29.20} are equivalent.

\end{remark}

This remark shows that without loosing generality
in the rest of the article
we impose

\begin{assumption}
                           \label{assumption 4.3.1}
We have $\rho_{b}=1$.
\end{assumption}

For    $\beta\geq 0$, introduce
Morrey's space $E_{p,q,\beta} $
as the set of $g\in  L_{p,q,\loc}$ such that  
\begin{equation}
                             \label{8.11.02}
\|g\|_{E_{p,q,\beta} }:=
\sup_{\rho\leq 1,C\in\bC_{\rho}}\rho^{\beta}
\dashnorm g  \|_{ L_{p,q}(C)} <\infty .
\end{equation}

Define
$$
E^{1,2}_{p,q,\beta} =\{u:u,Du,D^{2}u,
\partial_{t}u\in E_{p,q,\beta} \}
$$
and provide $E^{1,2}_{p,q,\beta} $ with an obvious norm.

It is important to have in mind that if
$\beta<2$ (our main case) and $u\in E^{1,2}_{p,q,\beta}$,
then according to Lemma 2.5 of \cite{Kr_23_1},
$u$ is bounded and continuous. 

Here is a useful approximation result.

\begin{lemma}[Lemma 2.3 of \cite{Kr_23_1}]
                               \label{lemma 3.14.3}
Let $f\in E_{p,q,\beta}$.
Define $f^{(\varepsilon)}$ as in the proof of Theorem 
\ref{theorem 10.21.1}.
Then for any $C\in \bC$ and $\beta'>\beta $
\begin{equation} 
                            \label{3.14.10}
\lim_{\varepsilon\downarrow0}
\|(f-f^{(\varepsilon)})I_{C}\|_{E_{p,q,\beta'} }=0.
\end{equation}
\end{lemma}

Introduce
$$
\cL u=\partial_{t}u+a^{ij} D_{ij}u +b^{i} D_{i}u .
$$

\begin{lemma}
                 \label{lemma 3.16.1}
 Let $x_{\cdot}$   be a solution  of
 \eqref{11.29.20} with $(t,x)=(0,0)$. Set
\begin{equation}
                                \label{3.16.5}
 p= p_{0}/\beta_{0}, \quad q= q_{0}/\beta_{0}
\end{equation}
and assume that

(a) for any $R\in(0,\infty)$
and Borel nonnegative $f$ on $\bR^{d+1}$
\begin{equation}
                             \label{12.11.06}
E\int_{0}^{\tau_{R}}f(s,x_{s})\,ds\leq N\|f\|_{E_{p,q,\beta_{0}}},
\end{equation}
where $N$ is independent of $f$ and $\tau_{R}$ is the first exit time of $(s,x_{s})$ from $C_{R}$. Then

(b)   for any $R\in(0,\infty)$ and   $u\in E^{1,2}_{p,q,\beta_{0}'} $, 
 with probability one for all $t\geq0$,
\begin{equation}
                                \label{3.16.2}
u(t\wedge\tau_{R},x_{t\wedge\tau_{R}})
=u( 0)+\int_{0}^{t\wedge\tau_{R}}D_{i}u 
\sigma^{ik} (s,x_{s})\,dw^{k}_{s}
+\int_{0}^{t\wedge\tau_{R}}\cL u (s,x_{s})\,ds
\end{equation}
and the stochastic integral above is a 
square-integrable
martingale.
\end{lemma}

Proof. By Corollary 5.6 of \cite{Kr_22} we have
$|Du|^{2}\in E_{r/2,s/2,2(\beta_{0}'-1) }$, where $r=
p\beta_{0}'/(\beta_{0}'-1)$,  $s=
q\beta_{0}'/(\beta_{0}'-1)$. Note that  
$$
  2>\beta_{0}'>1,\quad \beta_{0}'/(\beta_{0}'-1) >2,
\quad 2(\beta_{0}'-1)  < \beta_{0},
\quad
r /2\geq p, \quad s /2\geq q.
$$
This implies that the last statement of the lemma 
follows from  \eqref{12.11.06}.

To prove \eqref{3.16.2}, as in the proof of
Theorem \ref{theorem 12.11.3} write
$$
u^{(\varepsilon)}(t\wedge\tau_{R },x_{t\wedge\tau_{R }})
=u^{(\varepsilon)}( 0)+\int_{0}^{t\wedge\tau_{R }}D_{i}u^{(\varepsilon)} 
\sigma^{ik} (s,x_{s})\,dw^{k}_{s}
$$
\begin{equation}
                                      \label{3.16.3}
+\int_{0}^{t\wedge\tau_{R }}[
\partial_{t}u^{(\varepsilon)}+ a^{ij} D_{ij}u^{(\varepsilon)}
+b^{i} D_{i}u^{(\varepsilon)} ](s,x_{s})\,ds.
\end{equation}
Since $2(\beta_{0}'-1) <\beta_{0}$, by Lemma 
\ref{lemma 3.14.3} we have
$|Du-Du^{(\varepsilon)}|^{2}I_{C_{R}}\to0$
in $E_{p,q,\beta_{0}}$ and the stochastic integral
will converge in the mean square sense
as $\varepsilon\downarrow 0$ to the one
in \eqref{3.16.2}.

Since $u$ is bounded and continuous (Lemma 2.5 of \cite{Kr_23_1}), we have
the convergence  of the terms without integrals.
Regarding the integrals only the term with $b$
needs to be addressed. 

Notice that, thanks to $p= p_{0}/\beta_{0},q=
q_{0}/\beta_{0}$ and \eqref{3.14.1}
$$
\frac{d}{p }
+\frac{2}{q }\geq \beta_{0} \quad(>1).
$$
Also for any $C\in\bC$
$$
\dashnorm b\|_{L_{\beta_{0} p,\beta_{0} q}(C)}
=\dashnorm b\|_{L_{p_{0},  q_{0}}(C)}.
$$
This along with \eqref{12.11.06}
and  Remark 5.8 of \cite{Kr_22} imply that
$$
\|I_{C}b|Du-Du^{(\varepsilon)}|\,\|_{E_{p,q,\beta_{0}}}
\leq N\hat b\|I_{C}(u- u^{(\varepsilon)})\|_{E^{1,2}_{p,q,\beta _{0}}},
$$
where $N$ is independent of $\varepsilon$ and, owing to the fact that $\beta_{0}>\beta_{0}'$ and Lemma
\ref{lemma 3.14.3}, the right-hand side tends
to zero as $\varepsilon\downarrow0$.
This proves the lemma.\qed

We need the following fact which is a consequence of Theorem
\ref{theorem 10.27.1}.

\begin{theorem}
                              \label{theorem 3.14.1}
 Under Assumptions
\ref{assumption 12.12.3} 
and \ref{assumption 3.14.1} (with any $\hat b$)  there exists
$$
 \lambda_{0}  = \lambda_{0}  
(d,\delta,p_{0},q_{0},\beta_{0},\beta'_{0},\rho_{a},\alpha )>0
$$
such that    for $\gamma=\beta_{0}$ and $\gamma=\beta'_{0}$, $p= p_{0}/\gamma,  q= q_{0}/\gamma$, for any $\lambda \geq  \lambda_{0} $, Borel
$c(t,x)$ such that $ |c|\leq 1$,  
and $f\in E_{p,q,\gamma}$ there exists a unique
solution $u\in E^{1,2}_{p,q,\gamma}$ of
\begin{equation} 
                            \label{3.14.3}
\cL u-(\lambda+c)u+f=0.
\end{equation}
Furthermore for any $u\in E^{1,2}_{p,q,\gamma}$ we have
\begin{equation} 
                            \label{3.14.4}
\| \lambda u,\sqrt\lambda  Du, D^{2}u,\partial_{t}u\|_{E_{p,q,\gamma}}
\leq N_{0}\|\cL u-(\lambda+c)u\|_{E_{p,q,\gamma}},
\end{equation}
where $N_{0}=N_{0}( d,\delta,p_{0},q_{0},\beta_{0},\beta'_{0},\rho_{a},\alpha)$.
\end{theorem}

Actually,
Theorem
\ref{theorem 10.27.1}, in which \eqref{3.14.2} is replaced with 
$\mathfrak{b}_{\rho_{b}}\leq \tilde b$,
treats only the case
of $\gamma=\beta_{0}$. However, its assumptions are
also satisfied if we replace $\beta_{0}$ with $\beta_{0}'$. Then its conclusion holds true with such
a replacement as well, but for Theorem \ref{theorem 3.14.1}
to hold we need to take the largest of the $\check\lambda_{0}$'s
and the $N $'s corresponding to $\gamma=\beta_{0}$
and $\gamma=\beta_{0}'$ in Theorem \ref{theorem 10.27.1}.

In the following theorem we, in particular, specify 
the constant
$N$ in the definition of 
$\hat b$.

\begin{theorem}[Unconditional and conditional  weak uniqueness]
           \label{theorem 12.12.3}
Under Assumptions
\ref{assumption 12.12.3} 
and \ref{assumption 3.14.1} (with, generally, unspecified $\hat b$).

(i)   
If 
condition \eqref{4.15.3} is satisfied,
then   all solutions of \eqref{11.29.20} with fixed $(t,x)$ (provided they exist)
have the same finite-dimensional distributions.

(ii) Generally, let $y_{\cdot}$ and $z_{\cdot}$ be two solutions of
 \eqref{11.29.20} with $(t,x)=(0,0)$ perhaps on different probability spaces. Assume that
for $x_{\cdot}=y_{\cdot}$ and $x_{\cdot}=z_{\cdot}$
either (a)   or (b) 
  of Lemma \ref{lemma 3.16.1}
holds. 

Then $x_{\cdot}$ and $y_{\cdot}$ have the same finite-dimensional distributions.

\end{theorem}

Proof. First we prove (ii). Since by Lemma
\ref{lemma 3.16.1} (a) implies (b), we only need to show that (b) implies weak uniqueness.

Take   bounded Borel $c ,f $ on $\bR^{d+1}$
such that $0\leq   c \leq 1$.
By Theorem \ref{theorem 3.14.1} with $\lambda= \lambda_{0}$ there is
a bounded function $u$ {\em defined uniquely by\/} $a,b,c,\lambda,f$,
such that $u\in E^{1,2}_{p_{0}/\beta'_{0},q/\beta'_{0},\beta'_{0}}\subset E^{1,2}_{p_{0}/\beta_{0},q/\beta_{0},\beta'_{0}}$ and \eqref{3.14.3} holds.
In light of (b)
by   It\^o's formula  
applied to
$$
u(t,x_{t})\exp\Big(-\lambda t-\int_{0}^{t}c(s,x_{s})\,ds\Big)
$$
for any finite $T$
we obtain
$$
u(0)=Eu(T\wedge\tau_{R},x_{T\wedge\tau_{R}})\exp\Big(-\lambda (T\wedge\tau_{R})-\int_{0}^{T\wedge\tau_{R}}c(s,x_{s})\,ds\Big)
$$
$$
+E\int_{0}^{T\wedge\tau_{R}}f(t,x_{t})
\exp\Big(-\lambda t-\int_{0}^{t}c(s,x_{s})\,ds\Big)\,dt,
$$
where $\tau_{R}$ is the first exit time of $(s,x_{s})$
from $C_{R}$.
We send here $T,R\to\infty$ taking into account
that $\lambda+c\geq \lambda_{0}>0$, $u,f$ are bounded and $\tau_{R}\to\infty$. Then we get that
$$
E\int_{0}^{\infty}f(t,x_{t})
\exp\Big(-\lambda t-\int_{0}^{t}c(s,x_{s})\,ds\Big)\,dt
$$
  is uniquely defined by $a,b,c,\lambda,f$ (since it equals $u(0)$). For $T>0$ and $f=(\lambda+c)I_{t<T}$ this shows that
$$
E\exp\Big(-\lambda T-\int_{0}^{T}c(s,x_{s})\,ds\Big)
$$
is uniquely defined by $a,b,c,\lambda,T$. The arbitrariness
of $c$ and $T$ certainly proves  assertion (ii).

To prove (i), we take any solution
of \eqref{11.29.20}, say with $(t,x)=(0,0)$. 
 Set
$$
\frac{1}{\gamma}:=\frac{d}{p_{0}}+\frac{1}{q_{0}}
$$
and use
Theorem 4.2 of \cite{Kr_20} that
greatly simplifies in our situation. With
  $(p_{0},q_{0})/\gamma$ in place of its $(p_{0},q_{0})$,
in its conditional form this theorem yields that, since
for each $t\geq 0$, $\rho\leq 1$
and $C\in\bC_{\rho}$ we have that $\tau_{C}$, defined
as the first exit time of $(t+s,x_{t+s} )$, $s\geq0$,
from $C$, is less than $\rho^{2}$ and
since
$$
\|b\|_{L_{p_{0}/\gamma ,q_{0}/\gamma }(C)}^{2p_{0}/(p_{0}-\gamma d)}
\leq \big[N(d)(\hat b/\rho)\rho^{\gamma d/p_{0}+2\gamma /q_{0} }
\big]^{2p_{0}/(p_{0}-\gamma d)}
$$
$$=N(d,p_{0},q_{0})
\hat b^{2p_{0}/(p_{0}-\gamma d)}\rho^{2}
=N(d,p_{0},q_{0})
\hat b^{2q_{0}/ \gamma  }\rho^{2},
$$
we have with   constants $N$ depending only 
on $d,\delta,p_{0},q_{0}$, that
$$
E\Big\{\int_{0}^{\tau_{C}}|b(t+s ,x_{t+s})|\,ds\mid
\cF_{t}\Big\}
$$
$$
\leq N 
\big(\rho^{2}+\|b\|_{L_{p_{0}\gamma ,q_{0}/\gamma}(C)}^{2p_{0}/(p_{0}-\gamma d)}\big)^{\gamma d/(2p_{0})}\|b\|_{L_{p_{0}/\gamma ,q_{0}/\gamma }(C)}
$$
\begin{equation}
                               \label{3.8.3}
\leq N \rho^{\gamma d/p_{0}}(1+\hat b^{2q_{0}/\gamma})^{\gamma d/(2p_{0})}\hat b \rho^{
\gamma d/p_{0}+2\gamma/q_{0}-1}\leq N(d,\delta,p_{0},q_{0})
\hat b \rho,
\end{equation}
where in the last inequality we used that $\hat b\leq 1$. Hence, the left-hand side is less than $m_{b}\rho$,
provided that $N(d,\delta,p_{0},q_{0})
\hat b\leq m_{b}$. In that case all results of \cite{Kr_pta}
are available, in particular, estimate \eqref{12.11.6}
(implying \eqref{12.11.06}) 
holds for our solution whenever $d_{0}/p+1/q\leq 1$,
which in our case is guaranteed by the condition
$d/p_{0}+1/q_{0}\leq 1$ implying that for some $\beta_{0}\in(1,2)$
we have
$d_{0}/p_{0}+1/q_{0}=1/\beta_{0}$, $p_{0},q_{0}>\beta_{0}$ .
This proves 
the theorem.  \qed

Next, we proceed to proving the existence
of weak solutions.   Assumptions
\ref{assumption 12.12.3}, \ref{assumption 3.14.1},
and  \ref{assumption 4.3.1} are supposed to hold throughout the rest of this section
and we define $p,q$ as in \eqref{3.16.5}.
We start by drawing
consequences from Theorem
\ref{theorem 3.14.1}.  

\begin{corollary}
                            \label{corollary 3.14.1}
Assume that $a,b$ are smooth and bounded. Take
$R\leq 1$, smooth $f$, and let $u$ be the classical solution of
\begin{equation} 
                            \label{3.14.5}
\cL u+f=0
\end{equation}
in $C_{R}$ with zero boundary condition
on $\partial' C_{R}$.   Then
\begin{equation} 
                            \label{3.14.6}
|u|\leq NR^{2-\beta_{0}}\|I_{C_{R}}f\|_{E_{p,q,\beta_{0}}},
\end{equation}
where $N$ depends only on $d,\delta,p_{0},q_{0},\beta_{0},\rho_{a},\alpha$.
\end{corollary}

Indeed, the case $R<1$ is reduced to $R=1$ by using
parabolic dilations. If $R=1$, 
the maximum principle allows us to concentrate on $f\geq0$ and also shows that
$u(t,x)e^{\lambda _{0}t}$  is smaller  in $C_{1}$
than the solution $v$ of
$$
\cL v-\lambda _{0} v+I_{C_{1}}fe^{\lambda _{0}t}=0
$$
in $\bR^{d+1}$. Since $\beta_{0}<2$
by embedding theorems we have on $C_{1}$
$$
u\leq v 
\leq N\|v\|_{E^{1,2}_{p,q,\beta_{0}} }
\leq N\|I_{C_{\rho_{b}}}f\|_{E_{p,q,\beta_{0}}}.
$$

Now we can specify the constant $N_{1} $
in the definition of $\hat b$. Observe that so far
the size of $\hat b$ played a role only
in assertion (i) of Theorem \ref{theorem 12.12.3}.
\begin{corollary}
                            \label{corollary 3.14.2}
Assume that $a,b$ are smooth and bounded and let
$(\sft_{s},x_{s})$ be the corresponding Markov
diffusion process. Then for any $(t,x)\in\bR^{d+1}$,
 $\rho\leq1$, $C\in\bC_{\rho}$,  and Borel $f\geq0$
\begin{equation} 
                            \label{3.14.7}
I(t,x):=E_{t,x}\int_{0}^{\tau_{C}}f(t,x_{t})\,dt\leq N\rho^{2-\beta_{0}}\|I_{C }f\|_{E_{p,q,\beta_{0}}},
\end{equation}
where $\tau_{C}$ is the first exit time of 
$(\sft_{s},x_{s})$ from $C$. In particular,
\begin{equation} 
                            \label{3.14.8}
 E_{t,x}\int_{0}^{\tau_{C}}|b(t,x_{t})|\,dt\leq  N_{1}\rho \hat b,
\end{equation}
and in both estimates  $N$ and $N_{1}$ depend  only on $d,\delta,p_{0},q_{0},\beta_{0},\rho_{a},\alpha $.
\end{corollary}

Indeed, if $f$ is smooth, by It\^o's formula,
$I$ coincides with the solution of \eqref{3.14.5}
in a shifted $C_{\rho}$ and \eqref{3.14.7}
follows from \eqref{3.14.6}. For bounded Borel
$f$ we use the notation $f^{(\varepsilon)}$ from
the proof of Theorem 
\ref{theorem 10.21.1} and observe that $f^{(\varepsilon)}\to f$
almost everywhere, and the corresponding left-hand sides of \eqref{3.14.7} converge because they
are expressed in terms of the Green's function
of $\cL$. As far as the right-hand sides are concerned,
observe that by Minkowski's inequality
$\|f^{(\varepsilon)}\|_{E_{p,q}}\leq \|f \|_{E_{p,q}}$
and this yields \eqref{3.14.7} with $\|f \|_{E_{p,q}}$
in place of $\|fI_{C} \|_{E_{p,q}}$. Plugging 
$fI_{C}$ in such relation  in place of $f$ leads to \eqref{3.14.7}
as is. The passage to arbitrary $f\geq0$ is achieved
by taking $f\wedge n$ and letting $n\to\infty$.

To prove \eqref{3.14.8}
observe that due to self-similar transformations we may
assume that $\rho=1$, in which case we use \eqref{3.14.7} and the fact that
for $r\leq 1$ and $C'\in\bC_{r}$
$$
 r^{\beta_{0}}\dashnorm I_{C}b\|_{L_{p,q}(C')}\leq  r^{\beta_{0}}\dashnorm  b\|_{L_{p_{0},q_{0}}(C')}\leq Nr\dashnorm  b\|_{L_{p_{0},q_{0}}(C')}\leq N\hat b.
$$

Once $N_{1}$ is specified, we have the following.  

\begin{corollary}
                          \label{corollary 3.14.6}
Suppose that $a,b$ are smooth and bounded and let
$(\sft_{s},x_{s})$ be the corresponding Markov
diffusion process.  Then 
$ \sfb_{ 1}\leq m_{b}$ and all results
from \cite{Kr_pta} are applicable.
\end{corollary}

\begin{corollary}
                          \label{corollary 3.14.7}
Suppose that $a,b$ are smooth and bounded and let
$(\sft_{s},x_{s})$ be the corresponding Markov
diffusion process. Then for any $(t,x)\in\bR^{d+1}$,
Borel $f\geq0$, $T\in(0,\infty)$, there exists
$N$ depending only on $d,\delta$, $p_{0}$, $q_{0}$, $\beta_{0}$, $\rho_{a}$, $\alpha$, $T$, such that
\begin{equation}
                            \label{3.14.9}
 E_{t,x}\int_{0}^{T}f(t,x_{t}) \,dt\leq N
\|f\|_{E_{p,q,\beta_{0}}}.
\end{equation}
\end{corollary}

The proof of this is almost identical to the proof
of  \eqref{3.14.7} when $\rho=1$.

Here is a counterpart of Theorem \ref{theorem 10.21.1} in which $a$  and $b$ are not supposed to be smooth.
\begin{theorem}
               \label{theorem 3.15.1} 
(i) There is a probability space 
$(\Omega ,\cF ,P )$,
a filtration of $\sigma$-fields $\cF _{s}\subset \cF $, $s\geq0$,
a process $w _{s}$, $s\geq0$, which is a $d$-dimensional Wiener process
relative to $\{\cF _{s}\}$, and an $\cF _{s}$-adapted
process $x_{s}$ such that 
 (a.s.) for all   $s\geq0$ equation \eqref{11.29.20} holds with $(t,x)=(0,0)$.

(ii) Furthermore, for any 
  nonnegative Borel  $f$ on $\bR^{d+1} $ and
 $ T\in(0,\infty)$   we have
\begin{equation}
                                    \label{3.15.3}
E   \int_{0}^{T}  
f(s,x _{s} )\,ds \leq N  
 \|  f\| _{E_{p ,q,\beta_{0} }} ,
\end{equation}
where $N$ is the constant from \eqref{3.14.9}.

\end{theorem} 

Proof. As in the proof of Theorem 
\ref{theorem 10.21.1},
approximate $\sigma,b$ by smooth
$\sigma^{(\varepsilon)},b^{(\varepsilon)}$
ad take the corresponding Markov processes   
$(\sft_{t},x^{\varepsilon}_{t})$.
We noted in Corollary \ref{corollary 3.14.6} that
all results of \cite{Kr_pta} are available
for $(\sft_{t},x^{\varepsilon}_{t})$.
In particular, by Corollary 3.10 of \cite{Kr_pta}
for any $\varepsilon,n>0$ and
  $r>s\geq 0$  we have \eqref{11.24.4}
where $N=N(n,  d,\delta)$. This implies that the  $P_{0,0}$-distributions of $x^{\varepsilon}_{\cdot}$
are precompact on $C([0,\infty),\bR^{d})$
and a subsequence   $\varepsilon=\varepsilon_{n}
\downarrow 0$ of them converges
to the distribution of a process $x_{\cdot}=x^{0}_{\cdot}$ defined
on a probability space (the coordinate process on $\Omega=C([0,\infty),\bR^{d})$ 
with cylindrical $\sigma$-field $\cF$ 
completed with respect to
 $P$, which is
 the limiting distribution of $x^{\varepsilon}_{\cdot}$). Furthermore, by Theorem 5.1   of \cite{Kr_pta} for any 
  nonnegative Borel  $f$ on $\bR^{d+1} $ and
 $\varepsilon,T\in(0,\infty)$   we have
\begin{equation}
                                    \label{3.15.2} 
E _{0,0} \int_{0}^{T}  
f(s,x^{\varepsilon}_{s} )\,ds \leq N(d,\delta, T )  
 \|  f\| _{L_{d+1}(\bR^{d+1})} ,
\end{equation}
which by continuity is extended to $\varepsilon=0$
for bounded continuos $f$ and then by the usual measure-theoretic argument for all Borel $f\geq0$. 
After that estimate \eqref{3.15.2} also shows
that for any bounded Borel $f$ with compact support
\begin{equation}
                                    \label{3.15.5} 
\lim_{\varepsilon\downarrow 0}E _{0,0} \int_{0}^{T}  
f(s,x^{\varepsilon}_{s} )\,ds =
E _{0,0} \int_{0}^{T}  
f(s,x^{0}_{s} )\,ds .
\end{equation}
 
Furthermore, one can pass to the limit in
\eqref{3.15.3} written for $x^{\varepsilon}_{s}$
in place of $x_{s}$ and see that it holds
for $x^{0}_{s}$ if $f$ is bounded and continuous.
The extension of \eqref{3.15.3} to all Borel nonnegative $f$
is standard
and this proves assertion (ii).

 Now  we prove that assertions (i)  holds for
$x_{\cdot}$.
Estimate \eqref{11.24.4} implies that for any finite $T$
$$
\lim_{c\to\infty}P(\sup_{s\leq T}|x^{0}_{s}|>c)=0,
$$
and estimate   \eqref{3.15.3}
shows that for any finite $c$
$$
E   \int_{0}^{T}  I_{|x^{0}_{s}|\leq c}
|b (s,x^{0}_{s} )|\,dt <\infty.
$$
Hence, with probability one
$$
\int_{0}^{T}  
|b(s,x^{0}_{s} )|\,dt<\infty.
$$

Next, for $0\leq t_{1}\leq...\leq t_{n}\leq t\leq s$, bounded continuos $\phi(x(1),...,x(n))$,
and smooth bounded $u(t,x)$ with compact support by It\^o's formula we have 
$$
E_{0,0}\phi(x^{\varepsilon}_{t_{1}},...,x^{\varepsilon}_{t_{n}})
\Big[u(s,x^{\varepsilon}_{s})-u(t,x^{\varepsilon}_{t})-
\int_{t}^{s}\cL^{\varepsilon}u(r,x^{\varepsilon}_{r})\,dr\Big]=0,
$$
where
$$
\cL^{\varepsilon}u=\partial_{t}u+a^{\varepsilon ij}D_{ij}u+
b^{\varepsilon i}D_{i}u,\quad a^{\varepsilon}=(1/2)(\sigma^{(\varepsilon)})^{2}.
$$

Using \eqref{3.15.3}, Lemma \ref{lemma 3.14.3}, and the fact that $u$ has compact support show   that
$$
\lim_{\varepsilon_{1}\downarrow0}
\lim_{\varepsilon \downarrow0}E_{0,0}\int_{t}^{s}\big| b^{\varepsilon}-b^{(\varepsilon_{1})}\big|(r,x^{\varepsilon}_{r})|Du(r,x^{\varepsilon}_{r})|
\,dr = 0,
$$
$$
\lim_{\varepsilon \downarrow0}E_{0,0}\int_{t}^{s}  b^{(\varepsilon_{1})i} (r,x^{\varepsilon}_{r}) D_{i}u(r,x^{\varepsilon}_{r})
\,dr = E \int_{t}^{s}  b^{(\varepsilon_{1})i} (r,x^{0}_{r}) D_{i}u(r,x^{0}_{r})
\,dr,
$$
$$
\lim_{\varepsilon_{1} \downarrow0}E
\int_{t}^{s}\big| b -b^{(\varepsilon_{1})}\big|(r,x^{0}_{r})|Du(r,x^{0}_{r})|
\,dr = 0.
$$
 After that we easily conclude that
$$
E \phi(x^{0}_{t_{1}},...,x^{0}_{t_{n}})
\Big[u(s,x^{0}_{s})-u(t,x^{0}_{t})-
\int_{t}^{s}\cL u(r,x^{0}_{r})\,dr\Big]=0.
$$
It follows that the process
$$
u(s,x^{0}_{s}) -
\int_{0}^{s}\cL u(r,x^{0}_{r})\,dr
$$
is a martingale with respect to the completion
 of $\sigma\{x^{0}_{t}: t\leq s\}$. Referring to a well-known result from Stochastic Analysis proves assertion (i).
The theorem is proved. \qed

\begin{remark}
                           \label{remark 10.31.1}
In   \cite{RZ_20} the weak uniqueness is proved in the class
of solutions admitting, as they call it, Krylov type estimate
when $\sigma$ is {\em constant\/} and we have
$p,q\in[1,\infty]$ such that 
\begin{equation}
                           \label{10.29.5}
\frac{d}{p}+\frac{2}{q}=1,\quad\Big(
\int_{\bR}\Big(\int_{\bR^{d}} |b|^{p}\,dx\Big)^{q/p}
\,dt\Big)^{1/q}<\infty 
\end{equation}
(the Ladyzhenskaya-Prodi-Serrin condition).  

Actually, $p=\infty$, $q=2$ is not allowed in \cite{RZ_20}, this case fits in \cite{Kr_pt1} where
weak existence and conditional weak uniqueness is obtained. In case $p=d, q=\infty$ the comparison of the results in
\cite{Kr_pt1} and    \cite{RZ_20} can be found
in \cite{Kr_pt1}.

If $p\in [d+1,\infty)$ and we use the norms in \eqref{3.27.3}, set $p_{0}=p$ and $q_{0}:=q/2$ ($p\geq q_{0}$). Then for any $\rho>0$, and $C\in \bC_{\rho}$, by H\"older's inequality we have
\begin{equation}
                                 \label{10.29.6}
\|b\|_{L_{p_{0},q_{0}}(C)}
\leq \rho^{2/q}\|b\|_{L_{p,q}(C)},\quad
\dashnorm b\|_{L_{p_{0},q_{0}}(C)}
\leq N(d)\rho^{-1}\|b\|_{L_{p,q}(C)}
\end{equation}
and the last norm tends to zero as $\rho
\downarrow 0$. In that case also $d/p_{0}+1/q_{0}= 1 $.
This shows that in Assumption
\ref{assumption 3.14.1} the inequality $\mathfrak{b}
\leq\hat b\wedge \tilde b$ is satisfied 
  on the account of choosing $\rho_{b}$ small enough. 

By Theorem
\ref{theorem 12.12.3} (i) we get unconditional weak uniqueness of weak solutions that exist  by
Theorem \ref{theorem 3.15.1} if, say $\sigma$ is constant (as in \cite{RZ_20} and \cite{Ki_23}).

In case $d< p<d+1$ ($q>2p$), set $p_{0}=q_{0}=p$.
Then by H\"older's inequality we again get
the second relation in \eqref{10.29.6},
but this time our result is the same (if $\sigma$
is constant) as in
\cite{RZ_20} and \cite{Ki_23}: weak existence and conditional weak
uniqueness.

Interestingly enough, in case $p<d+1$ the estimates \eqref{10.29.6}
are still valid and show that  
\eqref{3.14.2} is satisfied for
an appropriate $\rho_{b}$, provided that the inequality in 
\eqref{10.29.5} is replaced with
\begin{equation}
                           \label{1.25.2}
\int_{\bR^{d}}\Big(\int_{\bR } |b|^{q}\,dt\Big)^{p/q}
\,dx<\infty.
\end{equation}
At the same time equations with $b$
satisfying \eqref{1.25.2} are not covered
in the literature so far (apart from the author's
works) and we   present in Remark \ref{remark 1.28.1} an example related to \eqref{1.25.2}.

\end{remark}

\begin{remark}
                  \label{remark 1.28.1}
There are examples showing that the assumptions
 of the second part of Theorem \ref{theorem 12.12.3}
concerning $b$ are satisfied with the norm in $L_{p,q}$
understood as in \eqref{1.25.2}
but not as in \eqref{10.29.5}
and \eqref{10.29.5}   does not hold no matter
what $p,q$ are, so that these examples are not
covered by  the results of
  \cite{RZ_20} or \cite{Kr_pt1}.  For instance,  
take $b(t,x)$ such that $|b|=cf$, where 
the constant $c>0$ and
$$
f(t,x)=I_{1>t>0,|x|<1}|x|^{-1}\Big(\frac{|x|}{\sqrt t}\Big)^{1/(d+1)}.
$$

If $p=\infty$, the second condition
in \eqref{10.29.5}, obviously, is not satisfied.
If $p=d,q=\infty$, so that the first
condition in \eqref{10.29.5} is satisfied, then
$$
\int_{|x|\leq 1}f^{d}(t,x)\,dx
=\int_{|x|\leq 1/\sqrt t}|x|^{-d+ d/(d+1)}\,dx
\to\infty
$$
as $t\downarrow 0$. Thus,   the second condition
in \eqref{10.29.5} is not satisfied
in this case. If $p\in(d,\infty)$ and $t>0$, $r<1$, then
$$
\int_{|x|\leq r}f^{p}(t,x)\,dx
=Nt^{(d-p)/2}\int_{0}^{r/\sqrt t}
\rho^{d-1-p+ p/(d+1)}\,d\rho=:
t^{(d-p)/2}I(r/\sqrt t).
$$
In order for that integral to converge,
we need 
\begin{equation}
                      \label{1.28.3}
p<d+1\quad (q>2(d+1)),
\end{equation}
and in this case $I(\rho)
\sim \rho ^{d-p+ p/(d+1)}$ as $\rho
\to\infty$.
Next,
$$
\int_{0}^{r^{2}}\Big(\int_{|x|\leq r}f^{p}(t,x)\,dx\Big)^{q/p}\,dt
=\int_{0}^{r^{2}}t^{(d-p)q/(2p)}
I^{q/p}(r/\sqrt t)\,dt
$$
$$
=2r^{2+(d-p)q/p}\int_{1}^{\infty}
\rho^{-3-(d-p)q/p}I^{q/p}(\rho)\,d\rho.
$$
Here the integrand has order of
$\rho^{-3+ q/(d+1)}$ and $-3+q/(d+1)
>-1$ by virtue of \eqref{1.28.3} and \eqref{10.29.5}. Therefore, the last integral diverges and condition 
\eqref{10.29.5} indeed fails to hold.

In contrast to this, although, if norms are understood
as in \eqref{3.27.3},   
the assumptions of  part (i) of Theorem \ref{theorem 12.12.3} concerning $b$,
obviously, are not satisfied   if $d/p_{0}+1/q_{0}=1$ and
$p_{0}\geq q_{0}$ (because then $p_{0}\geq d+1$), it turns out that they are satisfied with norms from \eqref{4.3.2} for some $q_{0}>p_{0}$ ($d/p_{0}+1/q_{0}=1,q_{0}>d+1$) if $c$ is small enough.

Indeed,  take $d+1<q_{0}<2(d+1)$
($p_{0}<d+1$) and note that for $r\le 3$
$$
\int_{0}^{r^{2}}f^{q_{0}}(t,x)\,dt
= I_{|x|<1}|x|^{-q_{0}}\int_{0}^{1\wedge r^{2}}
\Big(\frac{|x| }{\sqrt t}\Big)^{  q_{0}/(d+1)}  \,dt
$$
$$
\leq  |x|^{2- q_{0} }\int_{0}^{  r^{2} /|x|^{2}}t^{- q_{0}/(2d+2)}  \,dt
=:|x|^{2- q_{0} }J( r^{2} /|x|^{2}) ,
$$
where $J(s)\sim s^{1-q_{0}/(2d+2)}$ as $s\to\infty$.
Next,
$$
\int_{|x|<r}\Big(\int_{0}^{r^{2}}f^{q_{0}}(t,x)\,dt\Big)^{p_{0}/q_{0}}\,dx
$$
$$
= N
\int_{0}^{r}\rho^{d-1+(2- q_{0} )p_{0}/q_{0}}J^{p_{0}/q_{0}}(r^{2}/\rho^{2})\,d\rho
$$
$$
=Nr^{d+(2- q_{0} )p_{0}/q_{0}}
\int_{0}^{1}s^{d-1+(2- q_{0} )p_{0}/q_{0}}J^{p_{0}/q_{0}}(s^{-2})\,ds.
$$
Here in the integral with respect to $s$
$$
d+(2-q_{0})p_{0}/q_{0}-2\big(1-q_{0}/(2d+2)\big)p_{0}/q_{0}=d-p_{0}d/(d+1),
$$
which is strictly greater than zero and the above integral with respect to $ds$
is  finite implying that
\begin{equation}
                                \label{1.23.1}
\dashnorm f\|_{L_{p_{0},q_{0}}(C_{r})}\leq Nr^{-1}.
\end{equation}

If $r\leq 1$, $t\in (-r^{2},r^{2})$ and $|x|\leq 2r$, then
$$
\dashnorm f\|_{L_{p_{0},q_{0}}(C_{r}(t,x))}\leq
\dashnorm f\|_{L_{p_{0},q_{0}}(C_{3r} )}\leq Nr^{-1}.
$$
In case $r\leq 1$, $t\not\in (-r^{2},r^{2})$ or $|x|> 2r$ the left-hand side above is zero.
It follows that for small $c$ our $b$ satisfies
\eqref{3.14.2} with $\rho_{b}=1$. Also   $d/p_{0}+1/q_{0}=1$. Therefore,
  due to Theorem \ref{theorem 10.21.1},
\eqref{11.29.20}
has a   solution starting from any point, say if $\sigma$ is constant,
and all solutions starting from the same point have the same finite-dimensional distributions by Theorem \ref{theorem 12.12.3} (i).  

On the other hand, the results of \cite{Ki_23} still guarantee that there is weak existence and {\em conditional\/} weak
uniqueness in this example if $\sigma$ is constant.
Recall that our $\sigma$ is not necessarily constant or even continuous.

\end{remark}

By changing the origin we can apply Theorem \ref{theorem 12.12.3}
to prove the solvability of \eqref{11.29.20}
with any initial data $(t,x)$ and get  solutions with the
properties as in Theorems \ref{theorem 3.15.1} (ii)  
  weakly unique by Theorem
\ref{theorem 12.12.3}.
For such a solution denote by $P_{t,x}$ the distribution of $(\sft_{s},x_{s}),s\geq0$, ($\sft_{s}=t+s$)
on the Borel $\sigma$-field $\cF$ of $\Omega=C([0,\infty),\bR^{d+1})$. For $\omega=(\sft_{\cdot},x_{\cdot})\in \Omega$
set $(\sft_{s},x_{s})(\omega)=(\sft_{s},x_{s})$.
Also set $\frN_{s}=\sigma \{(\sft_{t},x_{t}),t\leq s\}$.

\begin{theorem}
                          \label{theorem 12.6.02}
The process
$$
X=\{(\sft_{\cdot},x_{\cdot}),\infty,\frN_{t},P_{t,x}\}
\quad (\infty\,\,\text{is the life time})
$$
is strong Markov with strong Feller resolvent
for which \eqref{9.27.5} holds true.
\end{theorem}

Proof. 
Take $u$ from Theorem \ref{theorem 3.14.1} with 
$\gamma=\beta_{0}'$,
 $c=0,
\lambda\geq\lambda_{0}$ and Borel bounded $f$. 
By It\^o's formula for any $(t,x) $
and $ 0\leq r\leq s$ we obtain that with
$P_{t,x}$-probability one  
$$
u(\sft_{s} ,x_{s})e^{-\lambda(s\wedge\tau_{R})} =u(\sft_{r },x_{r })e^{-\lambda(r\wedge\tau_{R})}
+\int_{r\wedge\tau_{R}}^{s\wedge\tau_{R}}e^{-\lambda v}
 \sigma^{ik}D_{i}u(\sft_{v},x_{v})
 \,dw^{k}_{v}
$$
\begin{equation}
                                    \label{12.14.02}
-\int_{r\wedge\tau_{R}}^{s\wedge\tau_{R}}
e^{-\lambda v}f(\sft_{v},x_{v})\,dv,
\end{equation}
where $\tau_{R}$ is the first exit time of 
$(\sft_{v},x_{v})$ from $C_{R}$

From \eqref{12.14.02} with $r=0$ as in the proof of 
Theorem \ref{theorem 12.12.3} we obtain
\begin{equation}
                                    \label{12.14.3}
E_{t,x}\int_{0}^{\infty}e^{-\lambda v}
f(\sft_{v},x_{v})\,dv=u(t,x).
\end{equation}
If $f$ is continuous, this implies that 
the Laplace transform of the continuous in $v$ function
$E_{t,x}f(\sft_{v},x_{v})$
is a  Borel function of $(t,x)$. Then the function
$E_{t,x}f(\sft_{v},x_{v})$ itself
is a  Borel function of $(t,x)$.  Since it is continuous
in $v$, it is Borel with respect to all its arguments.
 This 
fact is obtained for
bounded continuous $f$, but by usual measure-theoretic
arguments  carries it over to all Borel bounded $f$.

Then take $0\leq r_{1}\leq...\leq r_{m}=r$
and continuous $f$
and a  bounded Borel  function $\zeta\big(x(1),...,x(m)\big)$
on $\bR^{md}$  and conclude from \eqref{12.14.02} that
$$
E_{t,x}\zeta(x_{r_{1} },...,x_{r_{m} })u(\sft_{r} ,x_{r})e^{-\lambda r}
$$
$$
=E_{t,x}\zeta(x_{r_{1} },...,x_{r_{m} })\int_{r}^{\infty}e^{-\lambda v}
f(\sft_{v},x_{v})\,dv.
$$
In light of \eqref{12.14.3} this means that
$$
\int_{r}^{\infty}E_{t,x}\zeta(x_{r_{1} },...,x_{r_{m} })
e^{-\lambda v}E_{\sft_{r} ,x_{r}}f(\sft_{v-r},x_{v-r})\,dv
$$
$$
=\int_{r}^{\infty}E_{t,x}\zeta(x_{r_{1} },...,x_{r_{m} })e^{-\lambda v}
f(\sft_{v},x_{v})\,dv.
$$
We have the equality of two Laplace's transforms 
of functions continuous in $v$. It follows that
for $v\geq r$
$$
E_{t,x}\zeta(x_{r_{1} },...,x_{r_{m} })
 E_{\sft_{r} ,x_{r}}f(\sft_{v-r},x_{v-r})
=E_{t,x}\zeta(x_{r_{1} },...,x_{r_{m} }) 
f(\sft_{v},x_{v}).
$$
Again a measure-theoretic argument shows that
this equality holds for any Borel bounded $f$ and then
  the arbitrariness of $\zeta$ yields the Markov property
of $X$.

To prove that it is strong Markov it suffices to 
observe that, owing to \eqref{12.14.3}
 its resolvent $R_{\lambda}$ is strong Feller, that is maps bounded Borel functions into bounded continuous ones.

To deal with \eqref{9.27.5}, take, for instance,
$(t,x)=(0,0)$ and approximate our (conditionally weakly unique) solution
as in the proof of Theorem \ref{theorem 10.21.1}
by $x^{\varepsilon}_{\cdot}$.
for  
$R\in(0,\infty),y\in\bR^{d}$, 
introduce 
the functional $\gamma_{y,R}(x_{\cdot})$ 
 on $C([0,\infty),\bR^{d})$
 as the first exit time of $(s,x_{s})$
from $C_{R}(0, y)$. As is easy to see,
$\gamma_{y,R}(x_{\cdot})$ is lower semi-continuous.
It follows that the same is true
for
$$
 \int_{0}^{\gamma_{y,R}(x_{\cdot})}f( r,x_{r})\,dt,
$$
as  long as
a bounded continuous $f(t,x)\geq0$.
 It follows that  
\begin{equation}
                                \label{11.26.1}
\nliminf_{n\to\infty} E_{0,0} \int_{0}^{\gamma_{y,R}(x^{\varepsilon_{m}}_{\cdot})}f( r,x^{\varepsilon_{m}}_{r})\,dt
\geq
E_{0,0} \int_{0}^{\gamma_{y,R}(x^{0}_{\cdot})}f( r,x^{0}_{r})\,dt.
\end{equation}
In light of \eqref{3.15.3}, 
inequality
\eqref{11.26.1} holds   for $f=|b |
  $. If $f=|b|$
and $R\leq \rho_{b} $, as we have said in the proofs
of Theorem \ref{theorem 12.12.3} (i) and as it follows
from \eqref{3.14.8},
the  left-hand side of \eqref{11.26.1}
is smaller that $m_{b}R$. But then  
$$
E_{0,0}\int_{0}^{ \tau_{R}(y)}
|b( s, x^{0}_{ s}) |\,ds  \leq m_{b}R,
$$
and this with the possibility to change the origin leads to  \eqref{9.27.5}.
The theorem is proved. \qed

\mysection{A result from \protect\cite{Kr_23_1}}
                            \label{section 3.10.1} 

 The content  of this section is independent of Sections
\ref{section 12.24.1} and \ref{section 3.7.1},
however, we borrow some notation from 
Section \ref{section 3.7.1}.

We have  $p ,q ,\beta_{0}   $ such that
\begin{equation}
                        \label{3.21.01}
p ,q ,\beta_{0}\in(1,\infty), \quad \beta_{0}\in (1,2),
\quad \frac{d}{p }+\frac{2}{q }\geq\beta_{0}.
\end{equation}

Fix some $\rho_{a}  \in(0,\infty)$.
Parameters $\theta$ and $\check b$    below will be specified later.
\begin{assumption}
                         
We have
\begin{equation}
                                 \label{6.3.1}
 a^{\sharp}_{x,\rho_{a}}:=\sup_{\substack{\rho\leq\rho_{a}\\C\in\bC_{\rho}}}\dashint_{C}|a (t,x)-a_{C}(t)|
\,dx dt \leq  \theta,
\end{equation}
where
$$
a_{C}(t)=\dashint_{C}a (t,x)\,dxds\quad (\text{note}\quad t\quad \text{and}\quad ds).
$$
\end{assumption}

\begin{assumption}
                      
 We have 
\begin{equation}
                                    \label{4.7.1}
\mathfrak{b}_{1}:=\sup_{r\leq 1}r
\sup_{C\in \bC_{r}} 
\dashnorm b \|_{L_{p\beta_{0},q\beta_{0}}(C)}\leq   \check b .
\end{equation}
\end{assumption}

Let us specify $\theta$ in \eqref{6.3.1}.   It is easy to choose
$\theta_{1}(d,\delta,p,q)$ introduced in Lemma 4.5 of \cite{Kr_23_1}, so that it is a decreasing
function of $d$, and we suppose it is done.
In the following $\alpha\in(0,1)$ is a free parameter.

\begin{assumption}
                    \label{assumption 3.25.1}
For   $r$ defined as the least number such that
$$
 r\geq  (d+2)/\alpha,\quad r\geq p,q 
$$
and $\Theta(\alpha)=\{(p',q'):
p\leq p'\leq r,q\leq q'\leq r\}$
Assumption \ref{assumption 12.12.3} is satisfied with 
$$
\theta=\inf_{\Theta(\alpha)} \theta_{1}(d+1,\delta,p',q') =: 
\theta (d,\delta,p,q,\alpha).
$$
\end{assumption}

The fact that this $\theta>0$ is noted in \cite{Kr_23_1}.

Here is the main result of \cite{Kr_23_1} adjusted
to our needs.  
The constant $\nu=\nu(d,\beta_{0},p,q)$
below
is taken from Remark 2.2 of \cite{Kr_23_1}
when $\rho_{b }=1$.
\begin{theorem}
                       \label{theorem 10.27.1}
 
Under the above assumptions there exist 
$$
 \check b  = \check b (d,\delta, p, 
q, \rho_{a}, \beta_{0},\alpha )\in(0,1 ],
\quad \check \lambda_{0}=\check \lambda_{0}
(d,\delta, p, 
q, \rho_{a}, \beta_{0},\alpha)>0
$$
 such that,
if \eqref{4.7.1} holds with this $ \check b$,  
then for  any
$\lambda\geq   \check \lambda_{0} $,
function $c(t,x)$ such that $ |c|\leq 1$  and $f\in E_{p,q,\beta_{0} }$ there exists
a unique $E^{1,2}_{p,q,\beta_{0} }$-solution $u$
of $\cL u-(c+\lambda) u=f$. Furthermore, there exists
a constant $N$ depending only on  $d$, $\delta$, $p$, $q$, $\rho_{a}$,  $\beta_{0}$,
$\alpha$,      
such that 
\begin{equation}
                                \label{10.28.01}
\|\partial_{t} u,D^{2}u, \sqrt\lambda Du, \lambda u\|_{E_{p,q,\beta_{0} }}
\leq   N \nu^{-1}\|f\|_{E_{p,q,\beta_{0} }}.    
\end{equation}
\end{theorem}

\begin{remark}
                     \label{remark 10.27.2}   
The unique solution $u$ from Theorem \ref{theorem 10.27.1}
possesses the following properties

a) obviously, $u\in W^{1,2}_{p,q,\loc}$;

b) for $\beta_{0}>1$, by Lemma 2.6 of \cite{Kr_23_1}, we have
$Du\in L_{r,s,\loc} $, where $(r,s)=  (\beta_{0}-1)^{-1}\beta_{0}(p ,q )$;

c)   $u$ is bounded and continuous according to
Lemma 2.5 of \cite{Kr_23_1}.

\end{remark}

\end{document}